\documentclass{article}
\usepackage{amsfonts}
\usepackage{epsfig}
\begin{document}
\thispagestyle{plain}

\begin{center}
\Large{Some examples of Mahler measures as multiple polylogarithms}

\medskip
\normalsize{Matilde N. Lal\'{\i}n
{\footnote[1]{E-mail address: mlalin@math.utexas.edu}}
}{\footnote[2]{Supported by Harrington Fellowship}}

\medskip
University of Texas at Austin. Department of Mathematics.
 1 University Station C1200.  Austin, TX 78712, USA
\end{center}

\begin{abstract}
The Mahler measures of certain polynomials of up to five variables are given in terms of multiple polylogarithms. Each formula is homogeneous and its weight coincides with the number of variables of the corresponding polynomial.
\end{abstract}



\newcommand{\Cset}{\mathbb{C}}
\newcommand{\F}{\mathbb{F}}
\newcommand{\Rset}{\mathbb{R}}
\newcommand{\Qset}{\mathbb{Q}}
\newcommand{\Nset}{\mathbb{N}}
\newcommand{\Zset}{\mathbb{Z}}
\newcommand{\PP}{\mathbb{P}}
\newcommand{\HH}{\mathbb{H}}
\newcommand{\MU}{\mathbb{\mu}}
\newcommand{\TT}{\mathbb{T}}
\newcommand{\Li}{\mathrm{Li}} 
\newcommand{\I}{\mathrm{I}} 
\newcommand{\re}{\mathop{\mathrm{Re}}} 
\newcommand{\im}{\mathop{\mathrm{Im}}} 
\newcommand{\ii}{\mathrm{i}} 
\newcommand{\Lf}{\mathrm{L}} 
\newcommand{\dd}{\mathrm{d}}
\newcommand{\e}{\mathrm{e}}
\newcommand{\mod}{\mathrm{mod}}
\newcommand{\binom}[2]{{#1\choose#2}} 

\newcommand{\pf}{\noindent {\bf PROOF.} \quad}
\newcommand{\Res}{\mathrm{Res}}
\newcommand{\qed}{$\Box$}
\newtheorem{thm}{Theorem}
\newtheorem{defn}[thm]{Definition}
\newtheorem{prop}[thm]{Proposition}
\newtheorem{cor}[thm]{Corollary}
\newtheorem{lem}[thm]{Lemma}
\newtheorem{conj}[thm]{Conjecture}
\newtheorem{rem}[thm]{Remark}
\newtheorem{ex}[thm]{Example}
\newtheorem{exs}[thm]{Examples}
\newtheorem{obs}[thm]{Observation}

\newcommand{\T}{\mathbb{T}} 

\section{Introduction}

Lately there has been some interest in finding explicit formulae for the Mahler measure of polynomials. The
(logarithmic) Mahler measuxre of a polynomial $P \in \Cset [x_1, \dots,
  x_n]$ is defined as
\[ m(P) = \frac{1}{(2 \pi \ii )^n}\int_{\T^n} \log |P(x_1, \dots, x_n)|\, \frac{\dd x_1}{x_1},\, \dots, \frac{\dd x_n}{x_n}  \]
where $\T^n = \{(z_1,\dots,z_n) \in \Cset^n |\, |z_1| = \dots = |z_n| =1\}$ is the $n$-torus.

For the one-variable case, Jensen's formula
\[\frac{1}{2\pi} \int_0^{2\pi} \log | \e^{\ii \theta} - \alpha| \, \dd\theta = \log^+|\alpha|\]
(where $\log^+x = \log x$ if $x \geq 1$ and zero otherwise),  provides a simple expression
of the Mahler measure as a function on the roots of the  polynomial:
\[\mbox{Given} \quad P(x) = a \prod_j (x -\alpha_j),\quad \mbox{then}\quad m(P) = \log |a| + \sum_j \log^+|\alpha_j|\]

The two-variable case is much more complicated. Several examples
with explicit formulae have been produced.  Boyd \cite{B1,B2} and Smyth \cite{S1} have computed several examples and expressed some of them in terms of special values of L-functions of quadratic characters. Also Boyd \cite{B2} and Rodriguez Villegas \cite {RV}  have obtained analogous results with L-functions of certain elliptic curves.
Further,  Boyd and Rodriguxez Villegas \cite{BRV}, Maillot \cite{M},  and Vandervelde \cite{V} have produced examples where the Mahler measure is expressed as combinations of dilogarithms.

There are only a few exampxles for three variables. Smyth \cite{S2}
related the measure of $a+bx^{-1}+cy+(a+bx+cy)z$ to combinations of
trilogarithms and dilogarithms. Vandervelde \cite{V} obtained the
measure of $1+x+ az(x-y)$ as combinations of trilogarithms.

In this paper, we express the Mahler measure of some particular cases of up to five variable polynomials as combinations of multiple polylogarithms.  More precisely, we determine the Mahler measure of
\[ (1+w_1) \dots (1+w_n) + a (1-w_1) \dots (1-w_n) y\]
in $\Cset [w_1, \dots, w_n, y]$ for  $n=0,1,2,3$. We will refer to these as examples of the first kind.

We also consider
\[(1+w_1) \dots (1+w_n)(1+x) + a (1-w_1) \dots (1-w_n) (y+z)\]
in $\Cset [w_1, \dots, w_n, x, y, z]$ for $n=0,1,2$ (examples of the second kind).

In addition to these, we use the same method starting from Maillot's example \cite{M},  in order to compute the
Mahler measure of
\[(1+w)(1+y) + (1-w)(x-y)\]

\section{Summary of the results for the case $a=1$}
In order to be concrete and for future reference, we summarize the
results obtained for the particular case of $a=1$:

\renewcommand{\arraystretch}{2}
\begin{tabular}{|c|c|}
\hline
$ \pi m((1+x)+(1-x)y)$  &  $2\, \Lf(\chi_{-4}, 2)$\\
$\pi^2 m((1+w)(1+x)+(1-w)(1-x)y)$ &  $  7\, \zeta(3)$\\
$\pi^3 m((1+v)(1+w)(1+x)+(1-v)(1-w)(1-x)y)$  &  $ 7\, \pi \zeta(3) + 4   \sum_{0 \leq j < k} \frac{(-1)^j}{(2j+1)^2 k^2}$\\
\hline
$\pi^2 m((1+x)+(y+z)) $  & $  \frac{7}{2} \zeta(3)$\\
$\pi^3 m((1+w)(1+x)+(1-w)(y+z)) $ &  $  2 \pi^2 \Lf(\chi_{-4}, 2) +  8\sum_{0 \leq j < k} \frac{(-1)^{j+k+1}}{(2j+1)^3 k}$\\
$\pi^4 m((1+v)(1+w)(1+x)+(1-v)(1-w)(y+z)) $  & $ 93 \, \zeta(5)$\\
\hline
$ \pi^2 m((1+w)(1+y)+(1-w)(x-y))$  & $\frac{7}{2} \zeta(3) + \frac{\pi^2}{2} \log 2 $\\
\hline
\end{tabular}

The third and fifth formulae can be also written as

\[ \pi^3 m((1+v)(1+w)(1+x)+(1-v)(1-w)(1-x)y) \]
\begin{equation} \label{eq:L1}
= 7\, \pi \zeta(3) + 16 ( \Lf (\chi_{-4},\chi_0;2,2) - \Lf ( \chi_{-4}, \chi_{-4}^2;2,2))
\end{equation}

\[
\pi^3 m((1+w)(1+x)+(1-w)(y+z)) \]
\begin{equation}\label{eq:L2}
= \frac{7}{2}\pi \zeta(3) - 16 \log 2\, \Lf(\chi_{-4}, 3) +  16 ( \Lf ( \chi_0, \chi_{-4}; 1,3) - \Lf(\chi_{-4}^2, \chi_{-4}; 1,3))
\end{equation}

Here $\chi_0$ is the principal character and $\chi_{-4}$ is the real odd character of conductor 4, i.e.
\renewcommand{\arraystretch}{1.3}
\[\chi_{-4}(n) = \left \{ \begin{array}{rl}
1 & \quad\mbox{if}\quad n \equiv 1 \, \mbox{mod}\, 4\\
-1 & \quad\mbox{if}\quad n \equiv -1 \, \mbox{mod}\, 4\\
0 & \quad\mbox{otherwise}
\end{array} \right .
\]

The L functions are defined by
\[ \Lf(\chi_1,\dots,\chi_m;n_1,\dots,n_m) := \sum_{0<k_1<k_2< \dots <
  k_m} \frac{\chi_1(k_1)\chi_2(k_2) \dots \chi_m(k_m)}{k_1^{n_1}
  k_2^{n_2} \dots k_m^{n_m}}\]
This series is absolutely convergent if $ \re (n_m) >1$ and $\re (n_i)
  \geq 1$ for $i<m$.

We would like to point out that all these formulae are new, except for
$m((1+x)+(1-x)y)$, proved in \cite{S1}, and $m(1+x+y+z)$.

\section{Idea of the procedure and some technical steps}

Let $P_\alpha \in \Cset [x_1,\dots,x_n]$, a polynomial where the coefficients
depend polynomially on a parameter $\alpha \in \Cset $.
We replace $\alpha$ by $\alpha\, \frac{1-w}{1+w}$. A polynomial $\tilde{P}_\alpha \in \Cset [x_1, \dots , x_n, w]$
is obtained. The Mahler measure of the new polynomial is a certain
integral of the Mahler measure of the former polynomial. More
precisely,

\begin{prop}\label{lema} Let $P_\alpha \in \Cset [x_1,\dots, x_n]$ as above, then,
\begin{equation} \label{eq:lema1}
m(\tilde{P}_\alpha) = \frac{1}{2 \pi \ii } \int_{\T^1} m \left(P_{\alpha\, \frac{1-w}{1+w}}\right) \, \frac{\dd w}{w}
\end{equation}
Moreover, if the Mahler measure of $P_\alpha$ depends only on
$|\alpha|$, then
\begin{equation} \label{eq:lema2}
m(\tilde{P}_\alpha)= \frac{2}{\pi} \int_0^\infty m( P_x) \, \frac {|\alpha|\, \dd x}{x^2+|\alpha|^2}
\end{equation}
\end{prop}

\pf
Equality (\ref{eq:lema1}) is a direct consequence of the
definition of Mahler measure. In order to prove equality (\ref{eq:lema2}), write $w = \e^{\ii
  \theta}$. Observe that as long as $w$ goes through the unit circle
in the complex plane, $\frac{1-w}{1+w}$ goes through the imaginary
axis $\ii \Rset$, indeed, $\frac{1-w}{1+w} = -\ii
\tan\left(\frac{\theta}{2}\right)$. The integral becomes,
\[m(\tilde{P}_\alpha) = \frac{1}{2 \pi} \int_0^{2\pi} m \left(P_{\left|\alpha \, \tan\left( \frac{\theta}{2}\right) \right| } \right)
\, \dd\theta = \frac{1}{\pi} \int_0^{\pi} m \left(P_{|\alpha| \, \tan\left( \frac{\theta}{2}\right) } \right) \, \dd\theta \]

Now make $x = |\alpha| \tan \left(\frac{\theta}{2}\right)$, then $\dd\theta = \frac{ 2\, |\alpha|\, \dd x}{x^2+ |\alpha|^2}$,
\[m(\tilde{P}_\alpha) = \frac{2}{\pi} \int_0^{\infty} m \left(P_x\right) \frac{ |\alpha|\, \dd x}{x^2+ |\alpha|^2}\]
\qed

Our idea is to integrate the Mahler measure of some polynomials in
order to get the Mahler measure of more complex polynomials. We will need the following:

\begin{prop} \label{prop2} Let $P_a$ with $a>0$ be a polynomial as
  before, (its Mahler measure depends only on $|\alpha|$) such that
\[m(P_a) = \left \{ \begin{array}{ll} F(a) & \quad \mbox{if}\,\, a\leq1\\
G(a) & \quad \mbox{if}\,\, a > 1\end{array} \right. \]
Then
\begin{equation}
m(\tilde{P}_a) = \frac{2}{\pi} \int_0^1 F(x) \frac{a\, \dd x}{x^2+a^2} + \frac{2}{\pi} \int_0^1 G\left(\frac{1}{x}\right) \frac{a\, \dd x}{a^2x^2+1}
\end{equation}
\end{prop}

\pf
The Proof is the same as for equation (\ref{eq:lema2}) in Proposition \ref{lema},
with an additional change of variables $x
\rightarrow \frac{1}{x} $ in the integral on the right.
\qed

\bigskip

Recall the definition for polylogarithms which can be found, for instance, in Goncharov's papers, \cite{G1,G}:

\begin{defn} Multiple polylogarithms are defined as the power series
\[\Li_{n_1,\dots,n_m}(x_1,\dots, x_m) := \sum_{0 < k_1 < k_2 < \dots <k_m} \frac{x_1^{k_1}x_2^{k_2}\dots x_m^{k_m}}{k_1^{n_1}k_2^{n_2}\dots k_m^{n_m}}\]
which are convergent for $|x_i| < 1$. The weight of a polylogarithm function is the number $w=n_1+ \dots + n_m$.
\end{defn}

\begin{defn} Hyperlogarithms are defined as the iterated integrals
\[\I_{n_1,\dots, n_m}(a_1:\dots:a_m:a_{m+1}) := \]
\[\int_0^{a_{m+1}} \underbrace{ \frac{\dd t}{t-a_1}\circ \frac{\dd t}{t} \circ \dots  \circ \frac{\dd t}{t}}_{n_1} \circ \underbrace{\frac{\dd t}{t-a_2} \circ \frac{\dd t}{t} \circ  \dots \circ \frac{\dd t}{t}}_{n_2} \circ \dots \circ \underbrace{\frac{\dd t}{t-a_m}\circ  \frac{\dd t}{t} \circ \dots \circ \frac{\dd t}{t}}_{n_m}\]

where $n_i$ are integers, $a_i$ are complex numbers, and
\[\int_0^{b_{k+1}} \frac{\dd t}{t-b_1} \circ \dots \circ \frac{\dd t}{t-b_k} = \int_{0 \leq t_1 \leq \dots \leq t_k   \leq b_{k+1}} \frac {\dd t_1}{t_1-b_1}\, \dots \, \frac{\dd t_k}{t_k - b_k}\]
\end{defn}

The value of the integral above only depends on the homotopy class of the path connecting $0$ and $a_{m+1}$ on $\Cset  \setminus \{ a_1, \dots, a_m\} $. To be concrete, when possible, we will integrate over the real line.

It is easy to see (for instance, in \cite{G}) that,
\begin{eqnarray*}
\I_{n_1,\dots, n_m}(a_1:\dots:a_m:a_{m+1}) & = & (-1)^m \Li_{n_1,\dots, n_m}\left(\frac{a_2}{a_1}, \frac{a_3}{a_2},\dots,\frac{a_m}{a_{m-1}}, \frac{a_{m+1}}{a_m}\right)\\\\
\Li_{n_1,\dots, n_m}(x_1, \dots, x_m) & = & (-1)^m \I_{n_1,\dots,  n_m}((x_1 \dots  x_m)^{-1}: \dots :x_m^{-1}: 1)
\end{eqnarray*}
which gives an analytic continuation to multiple polylogarithms. For instance, with the above convention about integrating over a real segment, simple polylogarithms have an analytic continuation to $\Cset \setminus [1, \infty)$.

There are modified versions of these functions which are analytic in larger sets, like the Bloch-Wigner dilogarithm,
\begin{equation} \label{eq:defiD}
D(z) := \im( \Li_2(z)) + \log |z| \arg(1-z) \qquad z\in \Cset \setminus [1, \infty)
\end{equation}
which can be extended as a real analytic function in $\Cset \setminus
\{0, 1\}$ and continuous in $\Cset$.

We will eventually use some properties of $D(z)$:
\begin{equation} \label{eq:prop1D}
D(\bar{z}) = -D(z) \quad ( \Rightarrow D|_{\Rset} \equiv 0)
\end{equation}
\begin{equation} \label{eq:prop2D}
-2 \int_0^\theta \log|2 \sin t| \dd t =  D(\e^{2 \ii \theta}) = \sum_{n=1}^\infty \frac{\sin(2n\theta)}{n^2}
\end{equation}
More about $D(z)$ can be found in \cite{Z}.

Often we will write polylogarithms evaluated in arguments of
modulo greater than 1, meaning an analytic continuation given by the
integral. Although the value of these multivalued functions may not be uniquely
defined, we will always get linear combinations of these functions
which are one-valued, since they  represent Mahler measures of certain polynomials.

Now recall equation (\ref{eq:lema2}).  If the Mahler measure of $P_\alpha$ is a linear combination of multiple polylogarithms, and if we write $\frac{|\alpha|}{x^2+|\alpha|^2} =
\frac{\ii }{2}\left( \frac{1}{x+ \ii |\alpha|} - \frac{1}{x- \ii |\alpha|}\right )$, then it is
likely that the Mahler measure of $\tilde{P}_\alpha$ will be also a linear combination of
multiple polylogarithms. This will be the basis of our work.

In order to express the results more clearly, we will establish some
notation.

\begin{defn} Let  \[G:=  \left<\, \sigma_1,\, \sigma_2, \, \tau \right>\quad ( \cong\, \Zset/2\Zset \oplus \Zset/2\Zset \oplus \Zset/2\Zset)\] an abelian group generated by the following actions in the set $(\Rset ^*)^2$:
\begin{eqnarray*}
\sigma_1 : (a,b) & \mapsto & (-a, b)\\\\
\sigma_2 : (a,b) & \mapsto & (a, -b)\\\\
\tau : (a,b) & \mapsto & \left(\frac{1}{a}, \frac{1}{b}\right)
\end{eqnarray*}

Also consider the following multiplicative character:
\[\chi : G \longrightarrow \{-1, 1\}\]
\[\chi(\sigma_1) = -1\quad \chi(\sigma_2) = \chi(\tau) =1\]
\end{defn}

\begin{defn} Given $(a,b) \in \Rset^2$, $a \not = 0$, define,
\[\log (a,b) := \log |a|\]
\end{defn}

\begin{defn} Let $a \in \Rset ^* $, $x, y \in \Cset$,

\begin{eqnarray*}
{\mathcal L}^a_r(x) & := & \Li_r(xa) - \Li_r(-xa)\\\\
{\mathcal L}^a_{r:1}(x) & := & \log|a| ( \Li_r(xa) - \Li_r(-xa))\\\\
{\mathcal L}^a_{r,s}(x,y) & := & \sum_{\sigma \in G} \chi(\sigma) \Li_{r,s}\left( (x,y) \circ \left(a, \frac{1}{a}\right)^\sigma\right)\\\\
{\mathcal L}^a_{r,s:1}(x,y) & := & \sum_{\sigma \in G} \chi(\sigma)  \log \left(a, \frac{1}{a}\right)^{\sigma}  \Li_{r,s}\left((x,y) \circ \left(a, \frac{1}{a}  \right)^\sigma  \right)
\end{eqnarray*}

where \[(x_1,y_1) \circ (x_2, y_2) = (x_1 x_2, y_1 y_2)\]
is the component-wise product.
\end{defn}

\begin{obs} Let $a \in \Rset^*$, $x, y \in \Cset$, then,
\[{\mathcal L}_{r,s}^a(x,y) = {\mathcal L}_{r,s}^a(x, -y) = -{\mathcal  L}_{r,s}^a(-x, y)\]
and analogously with ${\mathcal L}^a_{r,s:1}$
\end{obs}

Observe also that the weight of any of the functions above is equal to the sum of its subindexes.

We will need some technical Propositions:

\begin{prop}\label{prop1} Given $a \in \Rset_{>0}$ 
\begin{equation}
\int_{\T^1} {\mathcal L}^{a\left|\frac{1-w}{1+w}\right|}_n(\ii )  \frac{\dd w}{w} = - 4 n {\mathcal L}^a_{n+1}(1) + 4 {\mathcal L}_{n:1}^a(1) 
\end{equation}
\end{prop}

\pf
By definition,
\begin{eqnarray}
\int_{\T^1} {\mathcal L}^{a\left|\frac{1-w}{1+w}\right|}_n( \ii ) \frac{\dd w}{w} & = &\int_{\T^1} \left(\Li_n\left( \ii a \left|\frac{1-w}{1+w}\right|\right) - \Li_n\left( - \ii a \left|\frac{1-w}{1+w}\right|\right)\right) \frac{\dd w}{w} \nonumber \\
& = & 2 \ii \int_0^\infty (\Li_n( \ii x) - \Li_n(- \ii x))\, \frac{2\,a\, \dd x}{x^2+a^2} \label{eq:prop1}
\end{eqnarray}
(which can be proved  in the same way as equation (\ref{eq:lema2}) in Proposition \ref{lema}).

Recall that 
\[\Li_n(\ii x) = - \I_n\left(\frac{1}{\ii x}:1\right) = - \int_0^1 \underbrace{\frac{\dd t}{t+\frac{\ii }{x}} \circ \frac{\dd t}{t} \circ \dots \circ \frac{\dd t}{t}}_n \]
using this and the fact that $\frac{1}{t+\frac{\ii }{x}} - \frac{1}{t-\frac{ \ii }{x}} = -\frac{2 \ii x}{t^2x^2 +1}$, the integral in (\ref{eq:prop1}) becomes
\begin{equation} \label{eq:2prop1}
-8 \int_0^\infty \int_{0 \leq t_1 \leq t_2 \dots \leq t_n \leq 1} \frac{x\, \dd t_1}{t_1^2 x^2+ 1} \, \frac{\dd t_2}{t_2}\, \dots \, \frac{\dd t_n}{t_n} \, \frac{a\, \dd x}{x^2+a^2}
\end{equation}

But 
\begin{eqnarray}
\int_0^\infty \frac{a\, x\, \dd x}{(t_1^2x^2 + 1)(x^2+a^2) } & = &  \int_0^\infty \left ( \frac{1}{x^2 + t_1^{-2}} - \frac{1}{x^2 + a^2} \right) \,\frac{a\, x\, \dd x}{t_1^2a^2-1} \nonumber \\
 & = & \frac{a}{2(t_1^2a^2 -1)} \, \left. \log\left( \frac{x^2+t_1^{-2}}{x^2+a^2}\right) \right | _0^\infty \nonumber \\
 & = & - \frac{a}{t_1^2a^2-1}\int_{t_1}^1 \frac{\dd s}{s} +\frac{a \log a}{t_1^2a^2-1} \label{eq:log}
\end{eqnarray} 
 
The integral in (\ref{eq:2prop1}) becomes
\[8 \int_{0 \leq t_1 \leq t_2 \dots \leq t_n \leq 1} \int_{t_1}^1 \frac{\dd s}{s} \frac{a\, \dd t_1}{t_1^2a^2-1}\,
 \frac{\dd t_2}{t_2}\, \dots \, \frac{\dd t_n}{t_n} -  8 \log a \int_{0 \leq t_1 \leq t_2 \dots \leq t_n \leq 1}  \frac{a\, \dd t_1}{t_1^2a^2-1}\, \frac{\dd t_2}{t_2}\, \dots \, \frac{\dd t_n}{t_n} \]
Although the sum of these two integrals is well defined, each of them
is not defined for $a >1$ if we choose the real segment $[0,1]$ as the
integration path. We choose a different integration path, such as $\gamma(\theta) = \frac{\e^{\ii \theta}+1}{2}$ for $-\pi \leq \theta \leq 0$.  

By using $\frac{2 a}{t_1^2a^2-1} = \frac{1}{t_1 - \frac{1}{a}} - \frac{1}{t_1 + \frac{1}{a}}$,  and $\int_{t_1}^1 \frac{\dd s }{s} = \int_{t_1}^{t_2} \frac{\dd s}{s} + \dots + \int_{t_n}^1 \frac{\dd s}{s}$
we get
\[4n \left(\I_{n+1}\left(\frac{1}{a}:1\right) - \I_{n+1}\left(-\frac{1}{a}:1\right) \right) -4 \log a  \left( \I_n\left(\frac{1}{a}:1\right) - \I_n\left(-\frac{1}{a}:1\right) \right)\]
\[ = -4n (\Li_{n+1}(a) - \Li_{n+1}(-a)) +4 \log a  (\Li_n (a) - \Li_n(-a)) \]
\[ =  -4n {\mathcal L}^a_{n+1}(1) +4 {\mathcal L}_{n:1}^a(1)\]
\qed

Observe that we will only use the above Proposition for the case $n=2$.

\begin{prop} \label{medium} For $z=\e^{\ii \theta}$, $a \in \Rset_{>0}$,
\begin{equation}
\int_0^1 {\mathcal L}^x_n(z) \frac{a}{x^2+a^2}\, \dd x + \int_0^1 {\mathcal L}^x_n(z) \frac{a}{a^2x^2+1}\, \dd x = \frac{\ii }{2} {\mathcal L}^a_{n,1}(\ii z, \ii )
\end{equation}
\end{prop}

\pf
By definition, the sum of the integrals is equal to
\begin{equation} \label{eq:medium1}
\int_0^1 (\Li_n(zx) -\Li_n(-zx)) \left( \frac{a}{x^2+a^2} + \frac{a}{a^2x^2+1} \right) \dd x
\end{equation}

Using that
\[\Li_n(zx) = - \int_0^1 \frac{\dd t}{t-\frac{1}{zx}} \circ \frac{\dd t}{t} \circ \dots \circ \frac{\dd t}{t} = - \int_{0 \leq t_1 \leq \dots \leq t_n \leq x} \frac{\dd t_1}{t_1 - \frac{1}{z}} \frac{\dd t_2}{t_2} \dots \frac{\dd t_n}{t_n}\]

The term in (\ref{eq:medium1}) with $\frac{a}{x^2+a^2}$ is equal to 
\[\int_0^1 \int_{0 \leq t_1 \leq \dots \leq t_n  \leq x} \left( \frac{1}{t_1 +\frac{1}{z}} - \frac{1}{t_1 -\frac{1}{z}}\right) \dd t_1 \, \,\frac{\dd t_2}{t_2}\, \dots \frac{\dd t_n}{t_n} \, \frac{a}{x^2 + a^2} \dd x \]

Writing $\frac{a}{x^2 + a^2} = \frac{\ii }{2} \left( \frac{1}{x+\ii a} - \frac {1}{x- \ii a}\right)$, we get 

\[\frac{\ii }{2}\left(\I_{n,1}\left(-\frac{1}{z}: -\ii a: 1\right) - \I_{n,1}\left(\frac{1}{z}: -\ii a: 1\right) +\I_{n,1}\left(\frac{1}{z}: \ii a: 1\right) -\I_{n,1}\left( -\frac{1}{z}: \ii a: 1\right) \right)\]
\[= \frac{\ii }{2}\left(\Li_{n,1}\left( \ii z a, \frac{\ii }{a}\right) - \Li_{n,1}\left(-\ii z a, \frac{ \ii }{a}\right) + \Li_{n,1}\left(\ii z a, -\frac{ \ii }{a}\right) - \Li_{n,1}\left(-\ii z a, -\frac{ \ii }{a}\right)\right) \]

The other integral can be computed in a similar way, (or taking advantage of the symmetry $a \leftrightarrow \frac{1}{a}$)
\[\frac{ \ii }{2}\left(\Li_{n,1}\left( \frac{\ii z}{a}, \ii a\right) - \Li_{n,1}\left( -\frac{\ii z}{a}, \ii a\right) + \Li_{n,1}\left(\frac{\ii z}{a}, -\ii a \right) - \Li_{n,1}\left(-\frac{\ii z}{a}, -\ii a \right)\right)\]

Adding both lines, we get the result.
\qed

\begin{prop}\label{complicated} For $z = \e^{\ii \theta}$, $a \in \Rset_{>0}$, 
\[\int_{\T^1} \int_0^1 {\mathcal L}^x_n(z) \frac{\left|a \frac{1-v}{1+v}\right|}{x^2+ \left|a \frac{1-v}{1+v}\right|^2}\, \dd x\, \frac{\dd v}{v} + \int_{\T^1} \int_0^1 {\mathcal L}^x_n(z) \frac{\left|a \frac{1-v}{1+v}\right|}{\left|a \frac{1-v}{1+v}\right|^2 x^2 + 1}\, \dd x\, \frac{\dd v}{v} \]
\begin{equation}
=2 \ii ({\mathcal L}^a_{n,2}(z,1)+ {\mathcal L}^a_{n,1:1}(z,1)) 
\end{equation}
\end{prop}

\pf
Consider the first integral. By definition, this is
\[ \int_{\T^1} \int_0^1 (\Li_n(zx)- \Li_n(-zx))\, \frac{\left|a \frac{1-v}{1+v}\right|}{x^2+ \left|a \frac{1-v}{1+v}\right|^2}\, \dd x\, \frac{\dd v}{v} \]
We do the same change of variables as in the Proof of equation (\ref{eq:lema2}) in Proposition \ref{lema} and we get
\begin{equation} \label{eq:complicated1}
 4 \ii \int_0^\infty \int_0^1 (\Li_n(zx) - \Li_n(-zx))\, \frac{y\, \dd x}{x^2+y^2}\, \frac{a\, \dd y}{y^2+a^2}
\end{equation} 

In the same way as in (\ref{eq:log}), we have:
\[\int_0^\infty \frac{a\,y\, \dd y}{(x^2+y^2)(y^2+a^2)} = \frac{a}{a^2 - x^2}\, \left( \log a + \int_x^1 \frac{\dd s}{s} \right)\] 
Integral (\ref{eq:complicated1}) becomes:
\[ - 4 \ii \int_0^1 \int_{0 \leq t_1 \leq \dots \leq t_n \leq x} \left( \frac{1}{t_1 + \frac{1}{z}} - \frac{1}{t_1 - \frac{1}{z}}\right) \dd t_1 \, \frac{\dd t_2}{t_2} \, \dots \, \frac{\dd t_n}{t_n}\, \left( \log a + \int_x^1 \frac{\dd s}{s} \right)  \frac{a\, \dd x}{x^2 - a^2}\]

This integral decomposes into two summands, one with $\int_x^1
\frac{\dd s}{s}$ and the other with $\log a$. But, as before,  when we do this, each
summand not longer converges if we integrate on the real interval
$[0,1]$ and if $0 \leq a \leq 1$. So, we will  change the path of
integration as we did before, to $\gamma(\theta) = \frac{\e^{\ii \theta} +1}{2} $ for $-\pi \leq \theta \leq 0$.

We first compute the integral with $\int_x^1 \frac{\dd s}{s}$. By using that $ \frac{a}{x^2-a^2} = \frac{1}{2} \left( \frac{1}{x-a} - \frac{1}{x+a}\right)$, we get
\[ 2 \ii \left(\I_{n,2}\left(\frac{1}{z}: a : 1 \right) - \I_{n,2}\left(-\frac{1}{z}: a : 1 \right) + \I_{n,2}\left(-\frac{1}{z}: -a : 1 \right) - \I_{n,2}\left(\frac{1}{z}: -a : 1 \right) \right)\]
\[ = 2 \ii \left(\Li_{n,2}\left(za , \frac{1}{a} \right) - \Li_{n,2}\left(-za , \frac{1}{a} \right) + \Li_{n,2}\left(za , -\frac{1}{a} \right) - \Li_{n,2}\left(-za , -\frac{1}{a} \right) \right)\]
\medskip

The term with $\log a$ yields
\[ 2 \ii \log a \left(\I_{n,1}\left(\frac{1}{z}: a : 1 \right) - \I_{n,1}\left(-\frac{1}{z}: a : 1 \right) + \I_{n,1}\left(-\frac{1}{z}: -a : 1 \right) - \I_{n,1}\left(\frac{1}{z}:-a : 1 \right) \right)\]
\[ = 2 \ii \log a \left(\Li_{n,1}\left(za , \frac{1}{a} \right) - \Li_{n,1}\left(-za , \frac{1}{a} \right) + \Li_{n,1}\left(za , -\frac{1}{a} \right) - \Li_{n,1}\left(-za , -\frac{1}{a} \right) \right)\]

The other integral is absolutely analogous, except that we use
\[\int_0^\infty \frac{a\, y\, \dd y}{(y^2x^2 + 1)(y^2 + a^2)} = \frac{a}{a^2x^2 -1}\left(\log a - \int_x^1 \frac{\dd s}{s} \right)\]
(we can also compute it using the symmetry $a \leftrightarrow \frac{1}{a}$).
\qed

\section{Examples of the first kind}
The Mahler measure of the polynomials that we study in this section
depends only on the absolute value of the parameter $\alpha$, hence, we will
write the formulae with $a = |\alpha|$ in order to simplify notation.

We start with the simple polynomial $1+ay$, whose Mahler measure is
\begin{equation} \label{eq:Jensen}
m(1+ay)= \frac{1}{2\pi} \int_0^{2\pi} \log |1+ a \e^{\ii \theta}| \, \dd\theta =  \log^+a
\end{equation}
The first application of our procedure yields:
\begin{thm} \label{easyteo} 
For  $a \in \Rset_{>0}$, 
\begin{equation}
\pi m((1+x) + a(1-x)y) =  -\ii {\mathcal L}^a_2( \ii )
\end{equation}
\end{thm}

\pf
By equation (\ref{eq:lema2}) in Proposition \ref{lema},
\begin{eqnarray}
 \pi m((1+x) +a(1-x)y) & = & \pi m(1+x) +  \pi m\left( 1 +a\frac{1-x}{1+x} y \right) \nonumber \\
 & = & 2 \int_0^\infty \log^+ z \, \frac{a\, \dd z }{ z^2+a^2} =  -\int_0^1 \log w \, \frac{2\,a\, \dd w}{w^2a^2 +1} \label{eq:thm12}\\
 & = & \ii \int_0^1 \int_w^1 \frac{\dd s}{s} \left ( \frac{1}{w+ \frac{\ii }{a} }- \frac{1}{w- \frac{\ii }{a}}\right) \dd w \nonumber 
\end{eqnarray}
 (we made $z= w^{-1}$).
\[= \ii \left( \I_2\left(-\frac{ \ii }{a}:1\right) - \I_2\left(
\frac{ \ii }{a}:1\right) \right) = - \ii \left( \Li_2\left( \ii a \right) - \Li_2\left(-\ii a\right)\right) = -\ii {\mathcal L}_2^a( \ii )\]
\qed

Recall that we mean the analytic continuation of $\Li_2$. If
we want to avoid this and work with the series, the formula should be
stated in the following way:

\begin{equation} \label{eq:split1}
\pi m((1+x) + a(1-x)y) =  \left \{ \begin{array}{ll} -\ii {\mathcal L}^a_2( \ii )  & \quad \mbox{if}\, a \leq 1\\ 
\pi \log a - \ii {\mathcal L}^{a^{-1}}_2( \ii )& \quad \mbox{if}\, a > 1 \end{array} \right.
\end{equation}

The $a \leq 1$ case is clear. For the $a > 1$ case,  
\begin{eqnarray*}
m( (1+x) + a(1-x)y) & = & m((1-x)y + a(1+x) ) \\
& = & \log a + m\left( \frac{1}{a}(1-x) y+ (1+x)\right)
\end{eqnarray*}
which proves the formula (\ref{eq:split1}).

\bigskip

Now we apply the procedure again:

\begin{thm} \label{thm13}
For $a \in \Rset_{>0}$, 
\begin{equation}
\pi^2 m((1+w)(1+x) + a(1-w)(1-x)y) = 4 {\mathcal L}^a_3(1) - 2 {\mathcal L}^a_{2:1}(1)
\end{equation}
\end{thm}

\pf
By Proposition \ref{prop1},
\begin{eqnarray*}
\pi^2 m((1+w)(1+x) + a(1-w)(1-x)y) & = & -\frac{1}{2}\int_{\T^1} {\mathcal L}_2^{a\left|\frac{1-w}{1+w}\right|}( \ii ) \frac{\dd w}{w}\\
& = & 4{\mathcal L}^a_3(1) - 2 {\mathcal L}^a_{2:1}(1)
\end{eqnarray*}
\qed

As before, we can express this with the following formula:

\[\pi^2 m((1+w)(1+x) + a(1-w)(1-x)y) = \left \{ \begin{array}{ll} 4 {\mathcal L}^a_3(1) - 2 {\mathcal L}^a_{2:1}(1)
  & \quad \mbox{if}\, a \leq 1\\  \pi^2 \log a + 4 {\mathcal L}^{a^{-1}}_3(1) - 2 {\mathcal L}^{a^{-1}}_{2:1}(1) & \quad \mbox{if}\, a > 1 \end{array} \right.\]

Note that we could compute the same Mahler measure using the formula
(\ref{eq:split1}) for the Mahler measure of 
$m((1+x) +a(1-x)y)$. By doing this, we obtain a different formula for
the Mahler measure of the polynomial considered in Theorem \ref{thm13}:

\begin{thm}\label{13}
For $a \in \Rset_{>0}$,
\[\pi^2 m((1+w)(1+x)+ a (1-w)(1-x)y) = -\ii \pi {\mathcal L}_2^a( \ii ) -{\mathcal L}^a_{2,1}(1, \ii )\]
\end{thm}

\pf
We integrate formula (\ref{eq:split1}) and use Proposition \ref{prop2}
\[ \pi^2 m((1+w)(1+x)+ a (1-w)(1-x)y)  =  - 2 \ii \int_0^1 {\mathcal L}_2^x( \ii ) \frac{a}{x^2 + a^2}\, \dd x \]
\begin{equation}
 + 2 \pi \int_0^1 \log\left(\frac{1}{x}\right) \frac{a}{a^2x^2 +1}\, \dd x - 2 \ii \int_0^1 {\mathcal L}_2^x(\ii ) \frac{a}{a^2x^2 + 1}\, \dd x \label{eq:eq1}
\end{equation}

The sum of the first and third integrals is ${\mathcal L}^a_{2,1}(-1, \ii ) = -{\mathcal L}^a_{2,1}(1, \ii )$ because of Proposition \ref{medium} with $n = 2$ and $z = \ii$.

The second integral is the same as the one that occurs in equation (\ref{eq:thm12}) and it yields $- \ii \pi{\mathcal L}^a_2(\ii )$.

Adding the three terms together we prove the statement. 
\qed

\bigskip

If we compare the two formulae that we have got for $m((1+w)(1+x)+ a
(1-w)(1-x)y)$, we have proved the following equality between multiple
polylogarithms:

\begin{cor} \label{coro15}
 For $a \in \Rset_{>0}$,
\begin{equation} \label{eq:coro}
4 {\mathcal L}^a_3(1) - 2 {\mathcal L}^a_{2:1}(1)= -\ii \pi{\mathcal L}^a_2(\ii ) - {\mathcal L}^a_{2,1}(1,\ii ) 
\end{equation}
\end{cor}

\pf
This Corollary is obtained from the two formulae for the Mahler
measure of $(1+w)(1+x)+a(1-w)(1-x)y$. See the Appendix for a direct proof.
\qed

\bigskip

Let us do the process of integration one more time.

\begin{thm} For $a \in \Rset_{>0}$,
\[ \pi^3 m((1+v)(1+w)(1+x) + a(1-v)(1-w)(1-x)y)\]
\begin{equation}
 = 4 \pi {\mathcal L}^a_3(1)- 2 \pi {\mathcal L}^a_{2:1}(1)- 2 \ii ({\mathcal L}^a_{2,2}( \ii, 1) + {\mathcal L}^a_{2,1:1}( \ii, 1))
\end{equation}
\end{thm}

\pf
We will integrate equation (\ref{eq:eq1}) of Theorem \ref{13}. The term of higher weight corresponds to 

\[-\int_{\T^1} \int_0^1 {\mathcal L}_2^x( \ii ) \frac{\left|a\, \frac{1-v}{1+v}\right|}{x^2 + \left|a\,\frac{1-v}{1+v}\right|^2} \dd x \frac{\dd v}{v}-\int_{\T^1} \int_0^1 {\mathcal L}_2^x( \ii ) \frac { \left| a \frac{1-v}{1+v}\right| }{\left|a \frac{1-v}{1+v}\right|^2 x^2 + 1}\,  \dd x\, \frac{\dd v}{v}\]
We solve this part by Proposition \ref{complicated}, setting $z = \ii$ and $n = 2$.

The third term is:
\[-\frac{\pi}{2} \int_{\T^1} {\mathcal L}_2^{a \left| \frac {1-v}{1+v}\right| }( \ii )  \frac{\dd v}{v} = 4 \pi {\mathcal L}_3^a(1)  - 2 \pi {\mathcal L}^a_{2:1}(1)\]
by Proposition \ref{prop1}.
\qed

\section{Examples of the second kind}
In this section we still have that all the Mahler measures only depend on $a=|\alpha|$.

We start with a different polynomial, $1+x+ay+az$. 

\begin{thm} 
For $a \in \Rset_{>0}$, 
\begin{equation} \label{eq:split3}
\pi^2 m(1+x+ay+az) = F(a) = \left \{ \begin{array}{ll} 2{\mathcal L}_3^a(1)  & \quad \mbox{if}\,\, a \leq 1\\  \pi^2 \log a + 2 {\mathcal L}_3^{a^{-1}}(1) & \quad \mbox{if}\,\, a \geq 1 \end{array} \right.
\end{equation}
\end{thm}

\pf
This was proved by Vandervelde \cite{V}. It is also possible to adapt some of the proofs of $m(1+x+y+z) = \frac{7}{2 \pi^2} \zeta(3)$.  For instance, following Smyth \cite{S2}, 
\[
\mbox{For} \quad V= 
\left(
\begin{array}{ccc}
2 & 0 & 1 \\
1 & 1 & 0 \\ 
1 & 1 & 1 
\end{array}
\right)
\qquad m(1+x+ay+az) = m((1+x+ay+az)^V)  \]
\begin{eqnarray*}
& = & m(1+x^2z+axy+axyz)  =  m(x^{-1} + ay + (x+ay) z) \\
& = & \frac{2}{\pi^2}( \Li_3(a) - \Li_3(-a)) =\frac {2}{\pi^2}{\mathcal L}_3^a(1)  
\end{eqnarray*}
for $0 \leq a \leq 1$. 

Another possibility is to adapt the elementary proof given in Boyd \cite{B1}. For $0 \leq a \leq 1$:
\[ \pi^2 m(1+x+ay +az) =  \pi^2 m(1+ay + x(1+aw)) = \pi^2 m\left( \frac{1+ay}{1+aw} + x \right)\]
 \begin{eqnarray*}
& = & \int_0^\pi \int_0^\pi \log^+ \left| \frac{1+a \e^{\ii t}}{1+a \e^{\ii s}}\right| \dd s \, \dd t =    \int_{0 \leq t \leq s \leq \pi} \log |{1+a \e^{\ii t}}| -\log |{1+a \e^{\ii s}}|\, \dd s \, \dd t  \\
& = &  \int_0^\pi (\pi-t)\, \log |{1+a \e^{\ii t}}| \dd t  -  \int_0^\pi s\, \log |{1+a \e^{\ii s}}| \dd s  = - 2 \int_0^\pi t\, \log |{1+a \e^{\ii t}}| \dd t  
\end{eqnarray*}
(here we have used that $0 \leq a \leq 1$, and formula (\ref{eq:Jensen})).

Now use that 
\begin{equation} \label{eq:sec}
 \log |1+a \e^{\ii t}| = \re \sum_{n=1}^\infty \frac{(-1)^{n-1}}{n}\,\, a^n \e^{\ii n t} = \sum_{n=1}^\infty \frac{(-1)^{n-1} \cos (n\,t)}{n} a^n 
\end{equation}
and apply integration by parts,
\[ \pi^2 m(1+x+ay+az)  =  - 2 \left. \sum_{n=1}^\infty \frac{(-1)^{n-1} \sin(n\,t)}{n^2}\,\, a^n t \,\right|^\pi_0 \]
\[ + 2 \int_0^\pi \sum_{n=1}^\infty \frac{(-1)^{n-1} \sin(n\,t)}{n^2}\,\, a^n \dd t =  4 \sum_{n =1 \, \mathrm{(odd)}}^\infty \frac{a^n}{n^3} = 2 (\Li_3(a) - \Li_3(-a)) \]

When $a \geq 1$, use that
\[m(1+x+ay+az)= \log a + m\left(\frac{1}{a} + \frac{x}{a} +y +z\right)\]
\qed

If we compare with formula (\ref{eq:split1}), for example, we may
wonder whether the formula in the second case $a \geq 1$ of
(\ref{eq:split3}) is a value of $2{\mathcal L}_3^a(1)$:

\[F(x) \stackrel{?}{=} 2{\mathcal L}_3^x(1),  \qquad x > 1\]

Meaning, as always, some branch of the analytic continuation of $\Li_3$.

We will see now that this is false. We should have for $x>1$, 
\[\pi^2 \log x + 2 \left(\Li_3\left(\frac{1}{x}\right) - \Li_3\left(-\frac{1}{x}\right)\right) \stackrel{?}{\equiv} 2 (\Li_3(x) - \Li_3(-x))\]

Differentiating, and using that $x\,\Li_3'(x) =\Li_2 (x)$
\[\frac{\pi^2}{x} + \frac{2}{x} \left (\Li_2\left(-\frac{1}{x}\right)  - \Li_2\left(\frac{1}{x}\right)  \right) \stackrel{?}{\equiv} \frac{2}{x} (\Li_2(x) - \Li_2(-x))\]

Multiplying by $x$, and differentiating again, 
\[\Li_1\left(\frac{1}{x}\right) - \Li_1\left( -\frac{1}{x}\right) = \Li_1(x) - \Li_1(-x)\]

Since $x>1$, the left term is 
\[\log \left(\frac{x+1}{x-1}\right)\]
(using that the principal branch of $\Li_1$ is equal to $- \log(1-x)$). The term in the right is equal to 
\[ - \int_0^1 \frac{1}{t-\frac{1}{x}} - \frac{1}{t+\frac{1}{x}} \, \dd t = -\lim_{\alpha \rightarrow \pi^-}\int_{-\alpha}^0
\frac{\frac{\ii \e^{\ii \theta}}{2}\, \dd\theta}{\frac{\e^{\ii \theta}+1}{2} - \frac{1}{x}}+\log(x+1) \] 
(we integrated on the path $\gamma(\theta) = \frac{\e^{\ii \theta}+1}{2}$ for $-\pi \leq \theta \leq 0$).
\[= -\log \left(\frac{x-1}{x} \right) + \lim_ {\beta \rightarrow \frac{\pi}{2}^-} \log \left(\frac{ \cos \beta\, \e^{-\ii \beta}\, x - 1}{x}\right) +\log(x+1) \]
\[= \log \left(\frac{x+1}{x-1}\right) - \ii \pi \]

$-\gamma(\theta) = \frac{\e^{- \ii \theta}+1}{2}$ for $-\pi \leq \theta \leq 0$ represents the other homotopy class and in this case, the integral is 
\[ = \log \left(\frac{x+1}{x-1}\right) + \ii \pi \]

Hence both functions are not equal, which implies that $F(x)$ can not
be expressed as $2{\mathcal L}_3^x(1) $.

\begin{thm}
For $a \in \Rset_{>0}$,
\begin{equation}
\pi^3 m((1+w)(1+x)+ a (1-w)(y+z))= -\ii \pi^2 {\mathcal L}^a_2( \ii )+ 2 \ii {\mathcal L}^a_{3,1}( \ii, \ii ) 
\end{equation}
\end{thm}

\pf 
Applying Proposition \ref{prop2} to formula (\ref{eq:split3}),
\[\pi^3 m((1+w)(1+x)+ a (1-w)(y+z)) = 4 \int_0^1{\mathcal L}_3^x(1)   \, \frac{a}{x^2+a^2}\, \dd x \]
\begin{equation} \label{eq:18}
+ 2\pi^2 \int_0^1 \log \left( \frac{1}{x}\right) \frac{a}{a^2x^2+1}\, \dd x + 4 \int_0^1 {\mathcal L}_3^x(1)  \frac{a}{a^2x^2+1}\, \dd x
\end{equation} 

The sum of the first and third integrals is $2 \ii {\mathcal L}^a_{3,1}( \ii, \ii )$ because of Proposition \ref{medium} with $n=3$ and $z=1$.

The second term is the same as the one in integral (\ref{eq:thm12}) in Theorem \ref{easyteo}, equal to $- \ii \pi^2{\mathcal L}^a_2( \ii )$.

Adding the three terms together we prove the statement.
\qed

\bigskip

Finally our result with the maximum number of variables:

\begin{thm} \label{thm:z5}
For $a \in \Rset_{>0}$,
\[ \pi^4 m((1+v)(1+w)(1+x) + a (1-v)(1-w)(y+z))  \]
\begin{equation}
 =4 \pi^2 {\mathcal L}_3^a(1)  - 2 \pi^2 {\mathcal L}^a_{2:1}(1) + 4\, ({\mathcal L}^a_{3,2}(1,1) + {\mathcal L}^a_{3,1:1}(1,1)) 
\end{equation}
\end{thm}

\pf
We will integrate equation (\ref{eq:18}) of the above Theorem as
always. The highest weight term corresponds to
\[ \frac{2}{ \ii } \int_{\T^1} \int_0^1 {\mathcal L}_3^x(1)  \frac{\left|a\, \frac{1-v}{1+v}\right|}{x^2 +\left|a\,\frac{1-v}{1+v}\right|^2} \dd x \frac{\dd v}{v} + \frac{2}{ \ii } \int_{\T^1} \int_0^1 {\mathcal L}_3^x(1)  \frac{\left|a\, \frac{1-v}{1+v}\right|}{\left|a\,\frac{1-v}{1+v}\right|^2x^2 + 1} \dd x \frac{\dd v}{v}\]
which can be evaluated using Proposition \ref{complicated} setting $z=1$, $n=3$.

Finally, 
\[ -\frac{\pi^2}{2} \int_{\T^1}{\mathcal L}_2^{a\left|\frac{1-v}{1+v}\right|}( \ii ) \frac{\dd v}{v} = 4 \pi^2 {\mathcal L}_3^a(1)  - 2 \pi^2 {\mathcal L}_{2:1}^a(1) \]
by Proposition \ref{prop1}.
\qed

\section{Integration of Maillot's formula}

So far we have been considering the cases of $1+ay$ and $(1+x) + a(y+z)$ and integrated them several times. We may wonder what happens with an "intermediate" case, namely $1+x+ay$. This case is not so easy to handle, so we will consider the following variant of the two variable case: $1+ax+(1-a)y$, with $a \in \Cset$. This time the Mahler measure will depend on the argument of $a$ as well. 

\medskip

According to Maillot \cite{M}, 
\begin{equation}
\pi m(a+bx+cy) = \left \{ \begin{array}{lr} \alpha \log |a| + \beta \log |b| +\gamma \log |c| + D\left(\left|\frac{a}{b}\right| \e^{\ii \gamma}\right) & \quad \triangle\\\\ \pi \log \max \{ |a|, |b|, |c|\} & \quad \mbox{not}\, \triangle \end{array} \right. 
\end{equation}

Where $\triangle$ stands for the statement that $|a|$, $|b|$, and
$|c|$ are the lengths of the sides of a triangle; and $\alpha$,
$\beta$, and $\gamma$ are the angles opposite to the sides of lengths
$|a|$, $|b|$ and $|c|$ respectively (Figure 1.A).

In our particular case,
\[ \pi  m( 1+ax+(1-a)y) = |\arg a| \log|1-a| + |\arg (1-a)| \log|a| \]
\begin{equation} \label{eq:Maillot}
+ \left \{ \begin{array}{cc} D(a)  & \quad \mbox{if} \, \im(a)\geq 0  \\   D(\bar{a}) & \quad \mbox{if} \, \im(a) <0
\end{array} \right . 
\end{equation}
since $|a|$, $|1-a|$ and 1 are the lengths of the sides of the
triangle whose vertices are $0$, $1$ and $a$ in the complex plane
(Figure 1.B). For the argument in the dilogarithm, $\gamma = |\arg(a)|$, then 
we have to take $a$ or $\bar{a}$ so $\gamma$ is always positive.

\begin{figure}
\centering
\epsfig{file=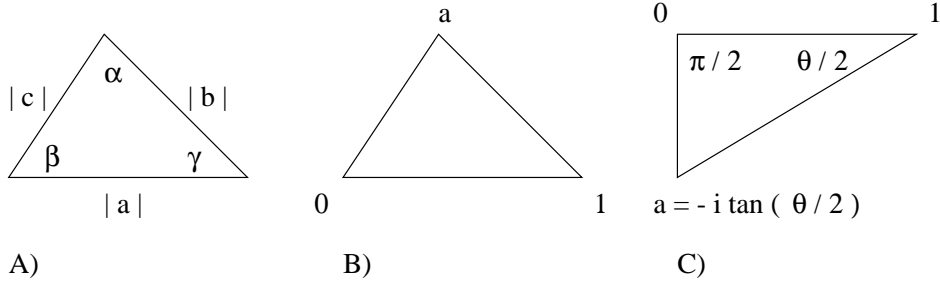}
\caption{A) Relation among the parameters in Maillot's formula.
B) Triangle for the general case of $1+a+(1-a)y$.
C) Triangle for the case $a = -\ii \tan \left( \frac{\theta}{2} \right)$.  }
\label{dibuM}
\end{figure}

We will integrate as always. We replace $a$ by $\frac{1-w}{1+w}$.

\begin{thm} \label{thm:M}
We have the following:
\begin{equation}
\pi^2 m((1+w)(1+y) + (1-w)(x-y))= 2{\mathcal L}^1_3(1) + \frac{\pi^2}{2} \log 2
\end{equation}
\end{thm}

\pf
We will apply Proposition \ref{lema} to equation (\ref{eq:Maillot}) and then
the change of variables $w = \e^{\ii \theta}$, which implies $ \frac{1-w}{1+w} = -\ii  \tan \left (
\frac{\theta}{2}\right)$. With that change, $a$ will be always pure
imaginary, so $|\arg a| = \frac{\pi}{2}$ (Figure 1.C).

\[ \pi^2 m((1+w)(1+y) + (1-w)(x-y))  = \]
\[ = \frac{1}{2} \int_{-\pi}^{\pi} \left( - \frac{\pi}{2} \log \left| \cos \left ( \frac{\theta}{2}\right) \right| + \left|\frac{\theta}{2}\right|  \log \left | \tan \left ( \frac{\theta}{2}\right) \right| + D\left(  \ii  \left| \tan \left ( \frac{\theta}{2}\right)\right| \right) \right) \, \dd\theta\]
Using the definition (\ref{eq:defiD}) of Bloch--Wigner dilogarithm,
\[= \frac{1}{2} \int_{-\pi}^{\pi} \left( - \frac{\pi}{2} \log \left| \cos \left ( \frac{\theta}{2}\right) \right| +\im \left(\Li_2 \left(\ii \left| \tan \left( \frac{\theta}{2}\right)\right| \right) \right) \right)\, \dd\theta  \]

For the first term, make $\tau = \frac{\pi - \theta}{2}$, this part becomes,
\begin{eqnarray*}
-\frac{\pi}{2}\int_0^\pi \log |\sin \tau| \dd\tau  & = & -\frac{\pi}{2}\int_0^\pi (\log |2 \sin \tau| - \log 2)  \dd\tau\\
& = & \frac{\pi}{4} D( \e^{2 \ii \pi}) + \frac{\pi^2}{2} \log 2 = \frac{\pi^2}{2} \log 2 
\end{eqnarray*}
by equation (\ref{eq:prop2D}).

For the second term, make $x = \tan \left( \frac{\theta}{2} \right)$, then $\dd\theta = \frac{2 \dd x}{x^2+1}$, the integral becomes,
\[ \int_{-\infty}^{\infty} \im(\Li_2( \ii |x|)) \frac{ \dd x}{x^2+1} = \frac{1}{\ii} \int_0^{\infty} (\Li_2(\ii x) - \Li_2(-\ii x)) \frac{\dd x}{x^2+1} = 2 {\mathcal L}^1_3(1)\]
The last equality is a particular case of the value computed for expression (\ref{eq:prop1}) in Proposition \ref{prop1}

Adding both terms we obtain the result.
\qed

\section{Concluding remarks}
To conclude, let us observe that all the presented formulae share a common feature. Let us assign weight 1 to any Mahler measure, to $\pi$ and to any logarithmic function.  Then all the formulae are homogeneous, meaning all the monomials have the same weight, and this weight  is equal to the number of variables of the corresponding polynomial.

\section*{Appendix 1: Some additional remarks about the
case $a=1$ }

In this appendix, we would like to give some additional details about the 
computation of the specific formulae that occur in the table of
results for $a=1$. All but two of these formulae can be directly deduced
from the Theorems we have proved. These exceptions are: the formula  with
the term $\zeta(5)$ and formula (\ref{eq:L2}).

For the formula with the term $\zeta(5)$:

\begin{thm} We have the following,
\begin{equation}\label{eq:z5}
\pi^4 m((1+v)(1+w)(1+x)+(1-v)(1-w)(y+z)) =  93 \, \zeta(5)
\end{equation}
\end{thm}

\pf By Theorem \ref{thm:z5} we have to prove that
\begin{equation} 
4 \pi^2 {\mathcal L}_3^1(1) + 4\,{\mathcal L}^1_{3,2}(1,1) = 93\,\zeta(5)
\end{equation}
i.e., 
\[
7 \pi^2 \zeta(3) + 8( \Li_{3,2}(1,1) - \Li_{3,2}(-1,1) +
\Li_{3,2}(1,-1) - \Li_{3,2}(-1,-1)) = 93\, \zeta(5)
\]
Now we use formula (75) of \cite{BBB}, which in this particular case,
states that
\[
\Li_{3,2}(x,y) = -\frac{1}{2} \Li_5(x\,y) + \Li_3(x)\, \Li_2(y) + 3
\Li_5(x) + 2 \Li_5(y) - \Li_2(x\,y) (\Li_3(x) + 2 \Li_3(y))  
\]
for $x, y = \pm 1$. 

Taking into account that
\begin{equation}\label{eq:z-1}
 \Li_k(1) = \zeta(k) \quad \mbox{and} \quad  \Li_k(-1) = \left( \frac{1}{2^{k-1}} - 1\right) \zeta(k)
\end{equation}
we get 
\[ \Li_{3,2}(1,1) - \Li_{3,2}(-1,1) +
\Li_{3,2}(1,-1) - \Li_{3,2}(-1,-1) =  - \frac{21}{4} \zeta(2) \zeta(3) + \frac{93}{8} \zeta(5)\]
We obtain the result by using that $\zeta(2) = \frac{\pi^2}{6}$
\qed

For formula (\ref{eq:L2}) we have the following 

\begin{prop} We have
\[ \sum_{0 \leq j < k} \frac{(-1)^{j+k+1}}{(2j+1)^3 k} = \frac{7}{4}
\zeta(3) \Lf(\chi_{-4},1) - \frac{3}{2} \zeta(2) \Lf(\chi_{-4},2) - 2
\log 2 \, \Lf(\chi_{-4},3) + 2 \sum_{0<m<n \,( m\; \mathrm{even)}}
\frac{\chi_{-4}(n)}{m \, n^3}\]
\end{prop}

\pf Writing $l=2j+1$ and $l+n=2k$, the left side is equal to
\[\sum_{0 \leq j < k} \frac{(-1)^{j+k+1}}{(2j+1)^3 k} = 2 \sum_{0 < l,
n \,  ( l, n\; \mathrm{odd})} \frac{
(-1)^{\frac{(l-1)+(l+n)}{2} +1 }}{l^3 (l+n)} = 2 \sum_{0 < l,
n \,  ( l\; \mathrm{odd)}} \frac{\chi_{-4}(n)}{l^3(l+n)}\]
Now write:
\[\frac{1}{l^3(l+n)} = \frac{1}{l^3\,n}- \frac{1}{l^2\,n^2} +
\frac{1}{l\, n^3}- \frac{1}{n^3(l+n)}\]
Using formulae (\ref{eq:z-1}), we get the first two terms of the
statement. We need to look at the last two terms together in order to ensure convergence,
\[ \sum_{0<n} \frac{\chi_{-4}(n)}{n^3} \sum_{0<l\,  ( l\;
\mathrm{odd})}  \left( \frac{1}{l} - \frac{1}{l+n}\right) = \sum_{0<n}
\frac{\chi_{-4}(n)}{n^3}  \left( - \log 2 + \sum^n_{0<m \,  ( m\;
\mathrm{even})} \frac{1}{m} \right)\]
and the statement follows.
\qed

Formula (\ref{eq:L2}) can be obtained from the above Proposition and
the fact that $\Lf(\chi_{-4},1) = \frac{\pi}{4}$.

\section*{Appendix 2: A direct proof for Corollary \ref{coro15}}

Recall the statement:


\begin{cor}
 For $a \in \Rset_{>0}$,
\begin{displaymath}
4 {\mathcal L}^a_3(1) - 2 {\mathcal L}^a_{2:1}(1)= -\ii \pi{\mathcal L}^a_2(\ii ) - {\mathcal L}^a_{2,1}(1,\ii )
\end{displaymath}
\end{cor}

\pf
First observe that after the change $a \leftrightarrow \frac{1}{a}$ the equality (\ref{eq:coro}) remains the same with the cancelation of a term of the form $\pi^2 \log a$ in each side. Then, it is enough to prove equation (\ref{eq:coro}) for $0 < a \leq 1 $.

Equation (\ref{eq:coro}) is equivalent to
\[ 8 \sum_{n = 1\, \mathrm{(odd)}}^\infty \frac{a^n}{n^3} - 4 \log a \sum_{n=1\, \mathrm{(odd)}}^\infty \frac{a^n}{n^2}  \stackrel{?}{=} - 2 \ii \pi \sum_{n=1\, \mathrm{(odd)}}^\infty \frac{\ii ^n\, a^n}{n^2}\]
\[+ \, 4 a \int_0^1 \frac{\dd t}{t^2+1}\, \circ \,\frac{\dd t}{t}\, \circ \, \left( \frac{\dd t}{t^2+a^2} + \frac{\dd t}{a^2t^2+1}\right)
\]

\[=  2 \pi \sum_{n=1\, \mathrm{(odd)}}^\infty \frac{ (-1)^{\frac{n-1}{2}} a^n}{n^2}\] 
\begin{equation} \label{eq:coro2}
 + 4 \int_0^1 \frac{\arctan s}{s} \left(\frac{\pi}{2} - \arctan \left( \frac{s}{a}\right) - \arctan (as) \right) \dd s 
 \end{equation}

Our strategy will be as follows: we will prove the equality for the derivatives and for the particular case $a=1$. In order to prove the equality for the derivatives, we will do the same, i.e., we will examine the case $a=1$ and differentiate again and compare the second derivatives.

Let us start with $a=1$. The term in (\ref{eq:coro2}) becomes

\[= 2 \pi \sum_{n =1 \, \mathrm{(odd)}}^\infty \frac{(-1)^{\frac{n-1}{2}}}{n^2} + 2\pi \int_0^1 \frac{ \arctan s}{s}\, \dd s - 8 \int^1_0 \frac {\arctan^2 s}{s} \, \dd s \]

Make $s = \tan x$:
\begin{eqnarray*}
\int^1_0 \frac {\arctan^2 s}{s} \, \dd s & = & \int_0^{\frac{\pi}{4}} \frac{x^2}{\sin x \cos x} \, \dd x \\
& = & \left . x^2\, \log ( \tan x) \right|_0^{\frac{\pi}{4}} - \int_0^{\frac{\pi}{4}} 2\,x\, (\log (2 \sin x) - \log(2 \cos x)) \dd x
\end{eqnarray*}
The first term is zero. For the second term make $y = \frac{\pi}{2} - x$
\begin{eqnarray*}
& = & - \int_0^{\frac{\pi}{4}} 2\,x\, \log (2 \sin x) \dd x+ \int_{\frac{\pi}{4}}^{\frac{\pi}{2}} (\pi - 2\,y) \,  \log(2 \sin y) \dd y \\
& = & - \int_0^{\frac{\pi}{2}} 2\,s\, \log (2 \sin s) \dd s + \pi \int_{\frac{\pi}{4}}^{\frac{\pi}{2}} \log (2 \sin s) \dd s
\end{eqnarray*}

Using properties (\ref{eq:prop1D}) and (\ref{eq:prop2D}), of the Bloch--Wigner dilogarithm,
\begin{eqnarray*}
& = & \left. s\, D(\e^{2 \ii s}) \right|_0^{\frac{\pi}{2}} - \int_0^{\frac{\pi}{2}} \sum _{n=1}^\infty \frac{\sin (2\,n\,s)}{n^2} \, \dd s - \frac{\pi}{2} \left . D(\e^{2 \pi \ii s }) \right|^{\frac{\pi}{2}}_{\frac{\pi}{4}}\\\
& = & - \sum_{n=1 \, \mathrm{(odd)}}^\infty \frac{1}{n^3} + \frac{\pi}{2} \sum_{n=1 \, \mathrm{(odd)}}^\infty \frac{(-1)^{\frac{n-1}{2}}}{n^2}
\end{eqnarray*}

Using the power series of $\frac{1}{1+s^2}$ and integrating, it is easy to see that
\[ \int_0^1 \frac{\arctan s}{s} \dd s= \sum_{n=1 \, \mathrm{(odd)}}^\infty \frac{(-1)^{\frac{n-1}{2}}}{n^2}\]

Putting all of this together in equation (\ref{eq:coro2}), we conclude:

\[-\ii \pi {\mathcal L}_2^1(\ii ) - {\mathcal L}^1_{2,1}(1, \ii ) =  8 \sum_{n=1 \, \mathrm{(odd)}}^\infty \frac{1}{n^3} = 4 {\mathcal L}^1_3(1)\]
as we expected.

We differentiate the original equation (\ref{eq:coro2}) and multiply it by $a$:
\[4 \sum_{n =1\, \mathrm{(odd)}}^\infty \frac{a^n}{n^2} - 4 \log a \sum_{n=1 \, \mathrm{(odd)}}^\infty \frac{a^n}{n} \]
\begin{equation} \label{eq:coro3}
\stackrel{?}{=} 2 \pi \sum_{n =1 \, \mathrm{(odd)}}^\infty \frac{\ii ^{n-1} a^n}{n} + 4 a \int_0^1 \arctan s \left( \frac{1}{s^2+ a^2} - \frac{1}{a^2s^2+1}\right) \dd s
\end{equation}

Set $a=1$, we get
\[4 \frac{3}{4} \zeta(2) - 2 \lim_{a \rightarrow 1} \log a ( \log (1+a) - \log (1-a)) \stackrel{?}{=}  \frac{\pi}{\ii} (\log(1+ \ii ) - \log(1- \ii ))\]
This is an equality, because the first term in the left and the term in the right are equal to $\frac{\pi^2}{2}$ and the other term is zero.

Apply integration by parts on the last term of (\ref{eq:coro3}):
\[4 a \int_0^1 \arctan s \left( \frac{1}{s^2+a^2} - \frac{1}{a^2s^2+1} \right) \dd s \]
\[ =  \pi \left( \arctan \left( \frac{1}{a}\right) - \arctan (a)\right) \\- 4 \int_0^1 \left( \arctan \left( \frac{s}{a} \right) - \arctan(as) \right) \frac{\dd s}{s^2+1} \]

Now we differentiate (\ref{eq:coro3}),
\[ -4 \log a \sum_{n =1 \, \mathrm{(odd)}}^\infty a^{n-1}\stackrel{?}{=} 2 \pi \sum_{n =1 \, \mathrm{(odd)}}^\infty ( \ii a)^{n-1}- \frac{2\, \pi}{a^2+1}\]
\[+ 4 \int_0^1 \left( \frac{1}{s^2+a^2} + \frac{1}{a^2s^2+1}\right) \frac{s\, \dd s}{s^2+1} \]
And this is an equality indeed, which can be seen from the fact:
\[\frac{4\, \log a}{a^2-1} = 4 \int_0^1 \left( \frac{1}{s^2+a^{-2}} - \frac{1}{s^2+ a^{2}} \right) \frac{s\, \dd s}{a^2-1}\]
\qed

\bigskip
\noindent{\bf Acknowledgements}

I am deeply grateful to Fernando Rodriguez Villegas for his constant
guidance and support and for sharing several ideas that have enriched
this work. I would also like to thank Sam Vandervelde for several
helpful discussions, and the Referee, whose suggestions led to several
improvements.


\begin{thebibliography}{9}
\bibitem{BBB} J. M. Borwein, D. M. Bradley, D. J. Broadhurst,
Evaluations of $k$-fold Euler/Zagier sums: a compendium of results for
arbitrary $k$, {\em Electronic J. Combin.\/} {\bf 4} (1997), no. 2,
\#R5.

\bibitem{B1} D. W. Boyd, Speculations concerning the range of Mahler's measure, {\em Canad.  Math. Bull.\/} {\bf 24} (1981), 453-469.

\bibitem{B2} D. W. Boyd, Mahler's measure and special values of L-functions, {\em Experiment. Math.\/} {\bf 7} (1998), 37-82.

\bibitem{BRV} D. W. Boyd, F. Rodriguez Villegas, Mahler's measure and the dilogarithm (I), {\em Canad. J. Math.\/} {\bf 54} (2002), 468-492.

\bibitem{G1} A. B. Goncharov, Polylogarithms in arithmetic and geometry, {\em Proc. ICM-94\/} Zurich (1995), 374-387.

\bibitem{G} A. B. Goncharov, Multiple polylogarithms and mixed Tate motives, 2001 (Preprint).

\bibitem{M} V. Maillot, G\'eom\'etrie d'Arakelov des vari\'et\'es toriques et fibr\'es en droites int\'egrables. {\em M\'em. Soc. Math. Fr. (N.S.)\/} {\bf 80} (2000) 129pp.

\bibitem{RV} F. Rodriguez Villegas, Modular Mahler measures I, in {\em Topics in number theory\/} (University Park, PA 1997), Math. Appl., 467, Kluwer Acad. Publ. Dordrecht, 1999, pp 17 - 48.

\bibitem{S1} C. J. Smyth, On measures of polynomials in several variables, {\em Bull. Austral. Math. Soc. Ser. A\/} {\bf 23} (1981), 49-63. Corrigendum (with G. Myerson): {\em Bull. Austral. Math. Soc.\/} {\bf 26} (1982), 317-319.

\bibitem{S2} C. J. Smyth, An explicit formula for the Mahler measure of a family of
3-variable polynomials, {\em J. Nombres Bordeaux\/} {\bf 14} (2002), 683-700.


\bibitem{V} S. Vandervelde, A formula for the Mahler measure of $axy+bx+cy+d$.  {\em J.
Number Theory\/} {\bf 100}  (2003),  no. 1, 184--202.


\bibitem{Z} D. Zagier, The Dilogarithm function in Geometry and Number Theory, {\em Number Theory and related topics, Tata  Inst. Fund. Res. Stud. Math.\/} {\bf 12} Bombay (1988), 231 - 249.
\end{thebibliography}
\end{document}